# THE PADÉ MATRIX PENCIL METHOD WITH SPURIOUS POLE INFORMATION ASSIMILATION

DANIEL TYLAVSKY[†], SONGYAN LI[†] and DI SHI[‡]

**Abstract.** We present a novel method for calculating Padé approximants that is capable of eliminating spurious poles placed at the point of development and of identifying and eliminating spurious poles created by precision limitations and/or noisy coefficients. Information contained in in the eliminated poles is assimilated producing a reduced order Padé approximant (PA). While the $[m + k/m]$ conformation produced by the algorithm is flexible, the $m$ value of the rational approximant produced by the algorithm reported here is determined by the number of spurious poles eliminated. Spurious poles due to coefficient noise/precision limitations are identified using an evidence-based filter parameter applied to the singular values of a matrix comprised of the series coefficients. The rational function poles are found directly by solving a generalized eigenvalue problem defined by a matrix pencil. Spurious poles place at the point of development, responsible in some algorithms for degeneracy, are identified by their magnitudes. Residues are found by solving an overdetermined linear matrix equation. The method is compared with the so-called Robust Padé Approximation (RPA) method [6] and shown to be competitive on the problems studied. By eliminating spurious poles, particularly in functions with branch points, such as those encountered solving the power-flow problem, solution of these complex-valued problems is made more reliable.

## 1   Nomenclature

$a_i(b_i)$—Numerator (denominator) coefficient
$d_i(d_i^{-1} = p_i)$—Inverse of value of pole $i$ (value of pole $i$)
c—vector of series coefficient
$c_i$—$i^{\text{th}}$ series coefficient
$e_i$—Numerator of the $i^{\text{th}}$ partial fraction
$l$–Pencil parameter
$m$–Denominator polynomial degree
$m + k$–Numerator polynomial degree
$n$—Number of series terms
$p$—Set of poles characterizing the function
$z = x + iy$—Complex argument variable.

## 2   Introduction.

When approximating a meromorphic function with a rational function, the singularities of the function appear as poles in the rational approximant, located near the singularities in the original function. The reverse, that all poles in the approximant represent singularities in the function, is not true. Poles not apparently related to the analytical properties of the function of interest are known as spurious poles and may accumulate in places where function convergence is expected, i.e., the convergence domain. An algorithm that eliminates all spurious poles while preserving the accuracy of the approximant would aid in convergence. However, the challenge of predicting of the location and occurrence of spurious poles remains elusive.

While the theory behind using exact arithmetic in the calculation of Padé approximants (arguably the most used and well understood method of obtaining a rational approximant), is well developed, predicting degeneracy and resolving numerical issues using finite precision remain [1]. Exact arithmetic can easily detect degeneracy, while finite precision arithmetic may mischaracterize degeneracy as ill-conditioning, leading the best-intentioned algorithms instead to solve the poorly condition problem, producing spurious poles.

---

[†] Dept. Electrical Computer and Energy Engineering, Arizona State University.
[‡] AINERGY LLC, Santa Clara, CA, USA.

Because of spurious poles, theoretical convergence for diagonal, $[m/m]$, Padé Approximants (PA's) is limited to convergence in capacity for function with branch points, which is much weaker than uniform convergence, though somewhat stronger than convergence in planar Lebesgue measure [2]-[4]. (Though convergence in capacity has been shown to imply point-wise convergence quasi everywhere for appropriate subsequences [5].) In any region absence of spurious pole, PA's converge uniformly. The goal of the present approach is to reduce the number of spurious poles, to zero in many cases, with minimal sacrifice to accuracy. By eliminating spurious poles, particularly in functions with branch points, such as those encountered solving the power-flow problem, solution of these complex-valued problems is made more reliable.

While the theory is scant on the origin of spurious poles, for the purposes here, we classify spurious poles into three categories:
- Poles place at the point of development.
- Poles accumulating on a disk centered at the radius of convergence.
- Poles appearing elsewhere in the convergence domain.

2.1    Spurious poles placed at the point of development

Assume the function, $f$, is holomorphic at the point of development, which unless stated otherwise is assumed to be the origin. (Similar definitions to the following may be established for functions developed at infinity.)

**Definition 1**: For each pair $k \in \mathbb{Z}, m \in \mathbb{N}$ with $(m + k) \in \mathbb{N}$, a rational function approximant is said to exist without defect if polynomials $P_{m+k,n} \in \mathbf{P}_{m+k}$ and $Q_{m+k,m} \in \mathbf{P}_m \setminus \{0\}$ exist such that.

$$(2.1) \qquad f(z) = \frac{P_{m+k,m}(z)}{Q_{m+k,m}(z)} + O(z^{2m+k+1}) \qquad z \to 0$$

where $O(\cdot)$ is Landau's big 'oh' and $\mathbf{P}_m$ the set of all complex polynomials of degree not greater than $m$.

By contrast a PA is defined as:

**Definition 2** [2]: For each pair $k \in \mathbb{Z}, m \in \mathbb{N}$ and $(m + k) \in \mathbb{N}$ if there exist two polynomials $p_{m+k,m} \in \mathbf{P}_{m+k}$ and $q_{m+k,m} \in \mathbf{P}_m \setminus \{0\}$ such that:

$$(2.2) \qquad q_{m+k,m}(z)f(z) - p_{m+k,m}(z) = O(z^{2m+k+1}) \qquad z \to 0$$

Then the rational function

$$(2.3) \qquad [m + k/m](z) := \frac{p_{m+k,m}(z)}{q_{m+k,m}(z)} \qquad z \to 0$$

is called the $[m + k, m]$-Padé approximant of the function $f$ (developed at zero) and is uniquely defined by (2.2).

Because the requirement of (2.2) results in $2m + k + 1$ homogeneous equations in $2m + k + 2$ variables, existence of a solution is guaranteed; however, because (2.1) is a nonlinear specification to the problem, it may not be possible to satisfy this equation for an arbitrary $k, m$. When the solution to (2.2) does not supply the accuracy requirement of $O(z^{2m+k+1})$ in (2.1), one or more poles are place at the point of development, the origin in this case. As we shall see, the Padé Matrix Pencil Method (PM$^2$) finds the poles of the rational-function approximant directly, allowing these poles to be easily screened out in the process.

A classic example of this behavior, which will be used to make an additional point in the next section, is the series in (2.4) along with its corresponding [1/1] PA [1], [6].

$$(2.4) \qquad f(z) = 1 + z^2 + O(z^3) \quad \to \quad [1/1] = \frac{z}{z} = 1$$

The order of the error in for this PA should be $O(z^3)$ but by inspection is $O(z^2)$. Clearly an $[m + k/m] = [1/1]$ PA cannot meet the error requirement of (2.1). The result is a pole and zero placed at the origin. Padé defined a deficiency (aka defect) index, $\gamma$, to measure the shortcomings of the $[m + k/m]$ selection, as the smallest integer for which,

(2.5) $$f(z) \neq \frac{p_{m+k,m}(z)}{q_{m+k,m}(z)} + O(z^{2m+k-\gamma_{k,m}+1}) \qquad z \to 0$$

For this example, $\gamma_{k,m} = 1$.

How we arrive at the solution to (2.2) matters a great deal if spurious poles are present. First we establish our notation for an arbitrary PA. Assuming a truncated series developed about zero of the form,

(2.6) $$f(z) = \sum_{i=0}^{2m+k} c_i z^i + O(z^{2m+k+1}) \qquad z \to 0$$

and the corresponding PA,

(2.7) $$f(z) = \frac{\sum_{i=0}^{m+k} a_i z^i}{\sum_{i=0}^{m} b_i z^i} + O(z^{2m+k+1}) \qquad z \to 0$$

the linear equations that define the PA for the case of interest $[m + k/m]$ are,

(2.8)
$$\begin{bmatrix} 0 \\ \vdots \\ 0 \end{bmatrix} = \begin{bmatrix} c_{m+k+1} & \cdots & c_{k+1} \\ \vdots & \ddots & \vdots \\ c_{2m+k} & \cdots & c_{m+k} \end{bmatrix} \begin{bmatrix} b_0 \\ \vdots \\ b_m \end{bmatrix} = Cb = 0 \qquad (a)$$

$$\begin{bmatrix} a_0 \\ a_1 \\ \vdots \\ a_j \\ \vdots \\ a_{m+k} \end{bmatrix} = \begin{bmatrix} c_0 & & & & \\ c_1 & c_0 & & & \\ \vdots & \vdots & \ddots & & \\ c_j & c_{j-1} & \cdots & c_0 & \\ \vdots & \vdots & & \ddots & \vdots \\ c_{m+k} & c_{m+k-1} & & \cdots & c_k \end{bmatrix} \begin{bmatrix} b_0 \\ \vdots \\ b_m \end{bmatrix} \qquad k \geq 0 \quad (b)$$

$$\begin{bmatrix} a_0 \\ a_1 \\ \vdots \\ a_{m+k} \end{bmatrix} = \begin{bmatrix} c_0 & & & \\ c_1 & c_0 & & \\ \vdots & \vdots & \ddots & \\ c_{m+k} & c_{m+k-1} & \cdots & c_0 \end{bmatrix} \begin{bmatrix} b_0 \\ \vdots \\ b_{m+k} \end{bmatrix} \qquad k < 0 \quad (c)$$

$$c_j = 0, j < 0 \qquad (d)$$

Once (2.8)(a) is solved, (2.8)(b)/(c) requires simply a vector-matrix multiplication. We refer to the algorithm that solves (2.8)(a) using the singular value decomposition approach as the SVD approach. If (2.8)(a) is solved using the matrix equation,

(2.9)
$$\begin{bmatrix} c_{m+k+1} \\ \vdots \\ c_{2m+k} \end{bmatrix} = \begin{bmatrix} c_{m+k} & \cdots & c_{k+1} \\ \vdots & \ddots & \vdots \\ c_{2m+k-1} & \cdots & c_{m+k} \end{bmatrix} \begin{bmatrix} b_1 \\ \vdots \\ b_m \end{bmatrix}$$
$$b_0 = 1$$
$$c_j = 0, j < 0$$

we refer to this as the Direct Method (DM) approach (cf. [1], Ch. 2).

While there are many methods of finding coefficients of (2.7), examining just two illustrates the requirements of an algorithm capable of eliminating spurious poles.

### 2.1.1 SVD Approach

Write (2.8)(a) succinctly as $Cb = 0$. Next we take the singular value decomposition (SVD) of $C$ in the form,

(2.10) $$C_{m \times m+1} = U_{m \times m} \Sigma_{m \times m+1} V^H_{m+1 \times m+1}$$

where $U$ and $V$ are unitary matrices, with superscript 'H' representing the conjugate transpose (Hermitian) operation and $\Sigma$ a diagonal matrix with real-valued entries. The last column (row) of $V$ ($V^H$) contains the solution vector, $b$.

When applied to calculate the $[1/1]$ of (2.4), our problem becomes,

(2.11) $$0 = \begin{bmatrix} c_2 & c_1 \end{bmatrix} \begin{bmatrix} b_0 \\ b_1 \end{bmatrix} = \begin{bmatrix} 1 & 0 \end{bmatrix} \begin{bmatrix} b_0 \\ b_1 \end{bmatrix}$$

The SVD of $C$ is given by,

(2.12) $$\begin{bmatrix} 1 & 0 \end{bmatrix} = [-1]\begin{bmatrix} 1 & 0 \end{bmatrix}\begin{bmatrix} -1 & 0 \\ 0 & 1 \end{bmatrix}$$

yielding values for $b_i$ from (2.8)(a) and values for $a_i$ consistent with (2.4) as shown in (2.13), which agrees with what we would get using the determinant method (cf. [1], Ch. 1.4).

(2.13) $$\begin{bmatrix} b_0 \\ b_1 \end{bmatrix} = \begin{bmatrix} 0 \\ 1 \end{bmatrix} \quad \begin{bmatrix} a_0 \\ a_1 \end{bmatrix} = \begin{bmatrix} 0 \\ 1 \end{bmatrix} \rightarrow [1/1] = \frac{z}{z}$$

This SVD approach result is consistent with the definition of spurious poles in [1], wherein an attempt is made to characterize them for the $[m/m]$ case. Spurious poles and zeros asymptotically (with increase number of terms) cancel each other out, as obviously happens in our example, though there are exceptions to this rule. Froissart doublets (discussed in the next major section) conform to this cancelation behavior [1].

### 2.1.2 Direct Matrix Approach

In the Direct Method (DM) approach, (cf. [1], Ch. 2), recognizing that the solution to (2.2) is only unique up to a scaling factor, the value of $b_0$ is taken as $b_0 = 1$. Solving (2.11) becomes solving the $1 \times 1$ matrix equation:

(2.14) $$[-c_2] = [c_1][b_1] \rightarrow [1] = [0][b_1]$$

where the coefficient matrix is singular for this degenerate case, leading to failure of the algorithm. Of course, it is easy to see that the assumption that $b_0 = 1$ has led to his degeneracy. Yet DM is one of the most widely used algorithms when the rational function form is needed as it is fast compared to the SVD approach.

### 2.1.3 Commentary

Degeneracies using DM can occur for nontrivial cases as well, using both exact and limited-precision arithmetic. For example, using exact arithmetic, all $[m/2]$ cases are degenerate for geometric series. Using limited precision on theoretically degenerate cases, the resultant roundoff error, which mimics noise added to the series coefficients (addressed in the next section), *may* mask the degeneracy, but at the expense of introducing spurious roots.

Using exact arithmetic and the SVD approach for degenerate cases will reveal one or more singular values of zero, which can be eliminated by reducing the degree of the numerator polynomial [6]. There are many other numerical approaches for finding PA's, and related approximants, (e.g., [12]), but many of these algorithms fail for these cases or because of numerical difficulties.

Like the SVD approach, as adapted in [6] for the elimination of spurious poles and zeros, PM$^2$ is able to avoid degeneracy, caused either by a short-sighted selection of $(k, m)$ or numerical degeneracy. We will look at the effects of noise and limited precision in the next section, but it is important to remember that the pole placed at the origin due to poor selection of $(k, m)$ is not a function of noise, since we observe this phenomenon using exact arithmetic, though noise can also place poles at the origin.

### 2.2 Spurious poles created by noise

Table 1 Poles and zeros differences using the SVD, DM and PM$^2$ approaches

|  | SVD | DM | PM$^2$ |
|---|---|---|---|
| Pole | $\varepsilon/(1 - \varepsilon^2/2) \approx \varepsilon$ | $\varepsilon$ | $\varepsilon$ |
| Zero | $\varepsilon/(1 - 3\varepsilon^2/2) \approx \varepsilon$ | $\varepsilon/(1 - \varepsilon^2) \approx \varepsilon$ | $\varepsilon/(1 - \varepsilon^2)$ |
| Separation | $\varepsilon^3$ | $\varepsilon^3$ | $\varepsilon^3$ |

Noise occurs in PA algorithms from multiple sources. The moment we move away from using exact arithmetic to using limited precision computing engines, truncation and roundoff error in the calculation of the series coefficients and then in the execution of the PA algorithm insert noise. In multi-point algorithms (see [1], Ch. 7) where measurements of the value of the function is the starting point, measurement noise (transducer, the A-to-D electronic conversion process including the resultant finite precision representation) along with practical uncertainty in the locations of the points in time and the location in the time-changing parameter-space add to the challenges of dealing with noise.

Consider again the function of (2.4) with noise of magnitude $\varepsilon$ added as shown in (2.15). When the SVD and the DM approaches are applied to the noisy form, this ceases to be a degenerate case for DM (provided $\varepsilon$ is not too small) and the spurious pole and zero at the origin have now been displaced, but the displacement are on the order of the inserted noise level, $\varepsilon$, as shown in Table 1.

$$(2.15) \qquad f(z) = 1 + \varepsilon z + z^2 + O(z^3)$$

There is limited theoretical work on characterizing spurious poles. Further, the fundamental causes of spurious poles are not well understood. The first work by Froissart [8] as reported and added to in [9], looked at the spurious poles for the geometric series, $\hat{g}(z)$, (2.16), where $g(z)$ was perturbed by noise $N(z)$, where $0 < \varepsilon \ll 1$ is the noise multiplier and where $r_i$ is a random draw. Where used numerically in [9], $r_i$ is taken as a uniform distribution on the interval $r_i \in [-1, 1]$.

$$(2.16) \qquad \begin{aligned} g(z) &= 1/(1-z) = \sum_{i=0}^{\infty} z^i \\ \hat{g}(z) &= g(z) + N(z) = \sum_{i=0}^{\infty} (1 + \varepsilon r_i) z^i \end{aligned}$$

This work revealed certain patterns in the characteristic pole placement and in the spurious pole locations [8], [9]. When studied for the PA's of conformation $[m + k/m], m, m + k \geq 0$, for a range of $k$ values, four types of roots were observed, as shown in Table 2.

i. One stable pole, $p_1$, whose distance from 1 was $|p_1 - 1| = O(\varepsilon)$ and decreased as $2m + k$ increased. This is a nonspurious pole that is characteristic of the underlying function.
ii. A cluster of (eponymous) Froissart doublets accumulating in pairs on the unit circle (radius of convergence) separated by $\varepsilon$.
iii. Poles or zeros, but not both simultaneously, at a distance $1/\varepsilon$ from the origin.

The Froissart doublets of (ii) are products of a synergy of noise and the choice to use a higher degree approximant than is supported by the data. These roots are paired by equal, but excess, degrees in the numerator and denominator.

The poles and zeros of (iii) are produced, again, by the synergy of noise and selection of an improper conformation of the PA, with the degree of discordance measured by its skewness from the latent conformation. More specifically, it is reported that any denominator degree in excess of the natural denominator-to-numerator degree ratio results in spurious poles, $p_i$, at a distance from the origin of $|p_i - 0| = O(1/\varepsilon)$ for $k \leq -1$. Likewise, any numerator degree in excess of the natural numerator-to-denominator degree ration results in spurious zeros, $z_i$, at a distance from the origin of $|z_i - 0| = O(1/\varepsilon)$ for $k > -1$.

All roots except the pole near 1 are spurious and ideally would be removed by a sufficiently sophisticated algorithm…without diminishing the accuracy of the PA significantly. Note that for this numerical experiment and for the previous one, the placement of the spurious poles is either related to the level of noise in the series signal or the poles accumulate on the circle in the complex plane whose radius is the radius of convergence (ROC).

Our duplication of these experiments using the DM algorithm and a 20-terms series shows that while the distance of the system pole from 1.0 is indeed $O(\varepsilon)$ (though this value can vary over an order of magnitude) and the number of Froissart doublets given by Table 2 is accurate, as is the number of Type (iii) poles/zeros, other metrics reported seem to be more aspirational. Froissart doublets can occur much further away from the unit circle than expected, their separation can be several orders of magnitude larger than expected and the Type (iii) poles/zeros in our simulations using the noisy geometric series occurred at a distance better described by $O(\log(\varepsilon^{-1}))$. (We used $\varepsilon \in [10^{-1}, 10^{-3}, 10^{-6}, 10^{-7}, 10^{-10}, 10^{-12}]$ in our simulations.)

The behavior of spurious poles becomes murkier as more complex functions and noise models are explored. For general rational functions of the type $[m - 1/m]$, with noise added in a way similar to (2.16), the placement of poles and zeros is not clear cut, even though some statements about the number of Froissart doublets remain applicable. Doublets still tend to accumulate on a circle, though the behavior is more chaotic [10]. The type of noise models used also impacts the predictability of the placement of spurious pole [11]. The statement that appears to hold in the most complex of situations is that, as $2m + k$ increases, the separation between Froissart doublets decreases, leading toward cancelation [2].

Table 2 PA roots for $\hat{g}(z)$, $m, m + k \geq 0$

| $k$ range | Poles at $\varepsilon$ from 1 | Poles at $1/\varepsilon$ | Zeros at $1/\varepsilon$ | Froissart Doublets near Unit Circle |
|---|---|---|---|---|
| $k < -1$ | 1 | $-k - 1$ | 0 | $(m + k)$ |
| $k \geq -1$ | 1 | 0 | $k + 1$ | $(m - 1)$ |

In all of these numerical experiments, because the perturbation of the series coefficients was much greater than the precision level of the computing engine used, the results shed some light on the effects of perturbations caused by precision limitations. Precision limitations, affecting the generation of series coefficients as well as roundoff in the numerics of the PA algorithm, can be modelled as noise in the series coefficients with calculations then proceeding in exact arithmetic. But the authors are unaware of any reliable models specific to the algorithms of interest. While amelioration of the effects of noise in the series coefficients may be possible for noise levels significantly greater than the precision level, ameliorating the effects of roundoff-level noise is more challenging. The goal of the present algorithm is to eliminate spurious poles due to noise, roundoff-error and unfortunate PA conformation choice while minimizing the impact of the accuracy of the PA.

### 2.3 Poles appearing elsewhere in the convergence domain

Of all the spurious pole phenomena, this is the least well understood phenomenon. Because these poles do not find their genesis in noise [2], but are sensitive to PA conformation as well as precision, these will likely remain the most difficult theoretically and numerically to address.

### 2.4 Organization

PM$^1$ is derived for the case $[m - 1/m]$ in the next section and its performance is demonstrated and compared to DM. The subsequent section shows the revised formulation, PM$^2$, in which filtering of spurious roots is applied, along with some short cuts that work in some cases. Then the algorithm is applied to the more general cases of $[m + k/m], m, m + k \geq 0$.

## 3  The Padé Matrix Pencil Method: Noiseless Data with Exact Arithmetic (PM$^1$)

Using PM$^1$, aka PM1, PA's with different conformations are handled differently. The derivation for PM$^1$ evolves most naturally for the $[m + k/m]$ case, with $k = -1$. This is discussed in the next subsection. The subsequent subsections address the cases for $k < -1$ and $k \geq 0$.

### 3.1  $k = -1$ Theory

For $k = -1$, we write the partial fraction expansion of the $[m - 1/m]$ PA as,

$$(3.1) \qquad f(z) = \frac{\sum_{i=0}^{m-1} a_i z^i}{\sum_{i=0}^{m} b_i z^i} + O(z^{2m}) = \sum_{j=1}^{m} \frac{e_j}{(1 - d_j z)} + O(z^{2m}) \quad z \to 0, \ 2m \leq n$$

where $d_j$ is the inverse of the $j^{\text{th}}$ pole, $p_j = d_j^{-1}$, of the partial fraction expansion with the corresponding residue of $-e_j/d_j$ and $n$ is the number of series coefficients. For this derivation, all poles are considered unique. (Non-uniqueness will be addressed later.) Now consider the denominator of the $j^{\text{th}}$ partial fraction.

$$(3.2) \qquad g_j(z) = (1 - d_j z)^{-1} \quad j \in \{1 \cdots m\}$$

The $i^{\text{th}}$ derivative is:

$$(3.3) \qquad g_j^{(i)}(z) = i! \, (d_k)^i (1 - d_k z)^{-i-1} \quad j \in \{1 \cdots m\}$$

The truncated power series expansion of $f(z)$,

(3.4)
$$f(z) \approx \sum_{i=0}^{n-1} c_i z^i = \sum_{i=0}^{n-1} \frac{f^{(i)}(z)}{i!} z^i \quad z \to 0$$

has the following coefficients:

(3.5)
$$c_i = \frac{f^{(i)}(0)}{i!} = \frac{1}{i!} \sum_{j=1}^{m} e_j \left((1 - d_j z)^{-1}\right)^{(i)}\bigg|_{z=0}$$
$$= \frac{1}{i!} \sum_{j=1}^{m} e_j (-1)^i (-d_j)^i i! (1 - d_j z)^{-i-1}\bigg|_{z=0} = \sum_{j=1}^{m} e_j d_j^i \quad i \in \{0 \cdots n - 1\}$$

It is easy though a bit tedious to show that the following is true:

(3.6) $$C_1 = D_1 E D_2$$
(3.7) $$C_2 = D_1 E D_0 D_2$$

where,

(3.8)
$$C_1 = \begin{bmatrix} c_0 & c_1 & \cdots & c_{l-1} \\ c_1 & c_2 & \cdots & c_l \\ \vdots & \vdots & \ddots & \vdots \\ c_{n-l-1} & c_{n-l} & \cdots & c_{n-2} \end{bmatrix}_{(n-l) \times l} \quad C_2 = \begin{bmatrix} c_1 & c_2 & \cdots & c_l \\ c_2 & c_3 & \cdots & c_{l+1} \\ \vdots & \vdots & \ddots & \vdots \\ c_{n-l} & c_{n-l+1} & \cdots & c_{n-1} \end{bmatrix}_{(n-l) \times l}$$

$$C_1^{ij} = c_{i+j}, \quad C_2^{ij} = c_{i+j+1},$$
$$i \in \{0, \cdots, n - l - 1\}, j \in \{0, \cdots, l - 1\}, m \le l \le n - m$$

(3.9)
$$D_1 = \begin{bmatrix} 1 & 1 & \cdots & 1 \\ d_1 & d_2 & \cdots & d_m \\ \vdots & \vdots & \ddots & \vdots \\ d_1^{n-l-1} & d_2^{n-l-1} & \cdots & d_m^{n-l-1} \end{bmatrix}_{(n-l) \times m}$$

(3.10)
$$D_2 = \begin{bmatrix} 1 & d_1 & \cdots & d_1^{l-1} \\ 1 & d_2 & \cdots & d_2^{l-1} \\ \vdots & \vdots & \ddots & \vdots \\ 1 & d_m & \cdots & d_m^{l-1} \end{bmatrix}_{m \times l}$$

(3.11) $$E = diag\{e_1, e_2, \cdots, e_m\}$$
(3.12) $$D_0 = diag\{d_1, d_2, \cdots, d_m\}, D_0^{-1} = diag\{p_1, p_2, \cdots, p_m\}$$

Consider the matrix pencil,

(3.13) $$C_1 - \lambda C_2 = D_1 E \{D_0^{-1} - \lambda I\} D_0 D_2$$

where $I$ is an $m \times m$ identity matrix and $D_0^{-1}$ is a diagonal matrix containing the poles of the PA. Observe that for any $\lambda = p_j, j = 1, \ldots, m$, the rank of $D_0^{-1} - \lambda I$ is diminished. Assuming the our value of $m$ leads to a PA of zero defect, it can be shown that, provided $m \le l \le n - m$, the poles of the PA may be found by solving the generalized eigenvalue problem [13],

(3.14) $$det(C_2^+ C_1 - \lambda I) = 0$$

where $C_2^+$ is the Moore-Penrose pseudoinverse of $C_2$,

(3.15) $$C_2^+ = [C_2^H C_2]^{-1} C_2^H$$

and where $C_2^H$ is the complex conjugate transpose of $C_2$. For our current specific case of interest, where $k = -1, n = 2m$, the dimension of the above matrices are all $m \times m \to n - l = l = m$, so we will not need the pseudo inverse; however, when spurious pole assimilation is used, we will need the flexibility of having $(n - l) > l \ne m$ and then the pseudoinverse will be needed.

Once the poles of the function are found, and assuming none of the poles are located at the origin (discussed later in PM[2]), the residues, $-e_j d_j^{-1}$, may be found using (3.5) in matrix form.

$$
(3.16) \quad \begin{bmatrix} c_0 \\ c_1 \\ \vdots \\ c_{n-1} \end{bmatrix} = \begin{bmatrix} 1 & 1 & \cdots & 1 \\ d_1 & d_2 & \cdots & d_m \\ \vdots & \vdots & \ddots & \vdots \\ d_1^{n-1} & d_2^{n-1} & \cdots & d_m^{n-1} \end{bmatrix}_{n \times m} \begin{bmatrix} e_1 \\ e_2 \\ \vdots \\ e_m \end{bmatrix} = c = De
$$

While we show the overdetermined form in (3.16), only the first $m$ equations need to be included if spurious pole assimilation (discussed later) is not used. Solving the overdetermined problem when unnecessary, while nonproblematic using exact arithmetic, limited precision may lead to unwanted precision-limitation-induced errors; further the $D$ matrix may become ill-conditioning and numerically singular if spurious pole assimilation is not used.

The formulation presented thus far does not account for poles with multiplicity greater than 1 and we have not encountered such cases in our application of the method. Without loss of generality, and to limit the complexity of the notation, consider the case of a function, $f(z)$, described by one pole, $d^{-1}$, of multiplicity, $p$.

$$
(3.17) \quad f(z) = \sum_{j=1}^{p} \frac{e_j}{(1 - dz)^j}
$$

The relationship between the $i^{\text{th}}$ Maclaurin series coefficient for $f(z)$ and the $j^{\text{th}}$ numerator value in (3.17) is given by,

$$
\begin{aligned}
c_i &= \frac{f^{(i)}(0)}{i!} = \frac{1}{i!} \sum_{j=1}^{p} e_j \left( (1-dz)^{-j} \right)^{(i)} \Big|_{z=0} \\
(3.18) \quad &= \frac{1}{i!} \left( e_1 d^i i! \, (1-dz)^{-i-1} \Big|_{z=0} + \cdots + e_j d^i \left( \frac{(i+j-1)!}{(j-1)!} \right) (1-dz)^{-i-1} \Big|_{z=0} \cdots \right) \\
&= \sum_{j=1}^{p} e_j d^i \left( \frac{(i+j-1)!}{i!\,(j-1!)} \right) \quad i = 0, \cdots, n
\end{aligned}
$$

For the case of $p = 3, m = 4$, the equivalent of (3.16) becomes:

$$
(3.19) \quad \begin{bmatrix} c_0 \\ c_1 \\ \vdots \\ c_7 \end{bmatrix} = \begin{bmatrix} 1 & 1 & 1 & 1 \\ d_1 & 2d_1 & 3d_1 & d_2 \\ \vdots & \vdots & \vdots & \vdots \\ d_1^7 & 8d_1^7 & 36d_1^7 & d_2^7 \end{bmatrix} \begin{bmatrix} e_1 \\ e_2 \\ e_3 \\ e_4 \end{bmatrix}
$$

In the subsequent equations developed, we shall ignore the possibility of poles with a multiplicity greater than 1, though the previous two equations may be inserted appropriately in our equations should pole multiplicity greater than 1 be of interest.

3.2    $k < 0$ Theory

The derivation of the previous section reveals the relationship between DM, PM[1] and the Gonnet et al. SVD approach, which then allows the general below diagonal and above diagonal cases to be arrived at more easily. If we reverse the order of the $b$ vector in (2.8)(a) and add a partition line for emphasis, we obtain the Hankel matrix:

$$
(3.20) \quad C\tilde{b} = 0 = \begin{bmatrix} c_{k+1} & c_{k+2} & \cdots & c_{k+m} & \vdots & c_{k+m+1} \\ c_{k+2} & c_{k+3} & \cdots & c_{k+m+1} & \vdots & c_{k+m+2} \\ \vdots & \vdots & \ddots & \vdots & \vdots & \vdots \\ c_{k+m} & c_{k+m+1} & \cdots & c_{2m+k-1} & \vdots & c_{2m+k} \end{bmatrix}_{m \times m+1} \begin{bmatrix} b_m \\ \vdots \\ b_0 \end{bmatrix}
$$

$$
c_j = 0, j < 0
$$

If we set $k = -1$, in (3.20), we see that everything to the left of the dashed line partitioning the $C$ matrix in (3.20), i.e., the first m columns, is the $C_1$ matrix in (3.8) and the last $m$ columns comprise the $C_2$ matrix in (3.8). The notation we use to represent this correspondence is:

(3.21) $$C = C_1 : C_2$$

This correspondence holds for arbitrary $k < 0$, and the algorithm proposed in the previous section may be implemented for $[m + k/m], -m \leq k < 0$ using the revise definitions for the $C_1$ and $C_2$ matrices given by (3.22). The equation for calculating the residues is unchanged from (3.16). For configuring the equations for this case, we again use $n - l = l, n = 2m + k + 1$ in (3.8)-(3.10).

(3.22)
$$C_1 = \begin{bmatrix} c_{k+1} & c_{k+2} & \cdots & c_{k+m} \\ c_{k+2} & c_{k+3} & \cdots & c_{k+m+1} \\ \vdots & \vdots & \ddots & \vdots \\ c_{k+m} & c_{k+m+1} & \cdots & c_{2m+k-1} \end{bmatrix}_{m \times m} \quad C_2 = \begin{bmatrix} c_{k+2} & c_{k+3} & \cdots & c_{k+m+1} \\ c_{k+3} & c_{k+4} & \cdots & c_{k+m+2} \\ \vdots & \vdots & \ddots & \vdots \\ c_{k+m+1} & c_{k+m+2} & \cdots & c_{2m+k} \end{bmatrix}_{m \times m}$$

$$[m + k/m], -m \leq k, n = 2m + k + 1$$

## 3.3 $k \geq 0$ Theory

Diagonal and above diagonal PAs are constructed as follows. For an $[m + k/m], k \geq 0, n = 2m + k + 1$ PA, the series is partitioned at term $c_{k+1}$, and $z^{k+1}$ factored out of the higher order terms.

(3.23)
$$f(z) \approx \sum_{i=0}^{k} c_i z^i + z^{k+1} \sum_{i=k+1}^{2m+k} c_i z^{i-k-1}$$

Using the techniques of the previous section, we then build an $[m - 1/m]$ for the higher-order term partition,

(3.24)
$$[m + k/m]_{f(z)} = \sum_{i=0}^{k} c_i z^i + z^{k+1} \frac{\sum_{i=0}^{m-1} \tilde{a}_i z^i}{\sum_{i=0}^{m} b_i z^i} \quad 2m + k + 1 = n, , k \geq 0$$

and by cross multiplying in (3.24) to achieve an $[m + k/m], k \geq 0$ PA.

(3.25)
$$[m + k/m]_{f(z)} = \frac{(\sum_{i=0}^{k} c_i z^i)(\sum_{i=0}^{m} b_i z^i)}{\sum_{i=0}^{m} b_i z^i} + z^{k+1} \frac{\sum_{i=0}^{m-1} \tilde{a}_i z^i}{\sum_{i=0}^{m} b_i z^i}$$

$$2m + k + 1 = n, , k \geq 0$$

The corresponding $C_1$ and $C_2$ matrices for calculating the denominator poles are the same as given by (3.22); however, as is obvious from (3.25), the corresponding residues are calculated for the fraction in (3.25) multiplied by $z^{k+1}$ and are given by:

(3.26)
$$\begin{bmatrix} c_{k+1} \\ c_{k+2} \\ \vdots \\ c_{n-1} \end{bmatrix} = \begin{bmatrix} 1 & 1 & \cdots & 1 \\ d_1 & d_2 & \cdots & d_m \\ \vdots & \vdots & \ddots & \vdots \\ d_1^{n-1} & d_2^{n-1} & \cdots & d_m^{n-1} \end{bmatrix}_{n \times m} \begin{bmatrix} e_1 \\ e_2 \\ \vdots \\ e_m \end{bmatrix} = c = De$$

It is clear from the development in section 3.1 that for the $[m - 1/m]$ case, provided the case is not degenerate, the denominator polynomial obtained is theoretically the same as that obtained by the matrix method, the difference being a scaling coefficient. That observation coupled with (3.22), indicates that the denominator polynomial obtained using (3.22) and the procedure of section 3.1 will match that of the matrix method, roundoff error notwithstanding. To show that the numerator polynomial matches that of the matrix method, we consider the case for $k \leq m$. The case for $k > m$ is addressed similarly.

Rather than solving (3.16) for the residues, $-e_j d_j^{-1}$, we may choose to solve for the numerator coefficients using the procedure of (2.8). The numerator coefficients $a_0 - a_k$ are found using the lowest order $k + 1$ terms of the cross-multiplication of the first term of (3.25), shown in matrix form immediately below.

(3.27)
$$\begin{bmatrix} a_0 \\ a_1 \\ \vdots \\ a_k \end{bmatrix} = \begin{bmatrix} c_0 & & & \\ c_1 & c_0 & & \\ \vdots & \vdots & \ddots & \\ c_k & c_{k-1} & \cdots & c_0 \end{bmatrix} \begin{bmatrix} b_0 \\ \vdots \\ b_k \end{bmatrix} \quad k \geq 0, k \leq m$$

The numerator coefficients for terms $a_{k+1} - a_{k+m}$ are found by solving for the numerator coefficients of the second term in (3.25) using (2.8) and by finishing the cross-multiplication of the first term, respectively.

$$(3.28) \quad \begin{bmatrix} \breve{a}_{k+1} \\ \breve{a}_{k+2} \\ \vdots \\ \breve{a}_{k+m} \end{bmatrix} = \begin{bmatrix} \tilde{a}_0 \\ \tilde{a}_1 \\ \vdots \\ \tilde{a}_{m-1} \end{bmatrix} = \begin{bmatrix} c_{k+1} & & & & \\ c_{k+2} & c_{k+1} & & & \\ \vdots & \vdots & \ddots & & \\ c_{k+m} & c_{k+m-1} & \cdots & c_{k+1} & 0 \end{bmatrix}_{m \times m+1} \begin{bmatrix} b_0 \\ \vdots \\ b_m \end{bmatrix}$$

$$k \geq 0, k \leq m$$

$$(3.29) \quad \begin{bmatrix} \hat{a}_{k+1} \\ \hat{a}_{k+2} \\ \vdots \\ \hat{a}_{k+m} \end{bmatrix} = \begin{bmatrix} 0 & c_k & c_{k-1} & \cdots & c_{k-m+1} \\ 0 & 0 & c_k & \cdots & c_{k-m+2} \\ 0 & \vdots & \vdots & \ddots & \\ 0 & 0 & 0 & \cdots & c_k \end{bmatrix}_{m \times m+1} \begin{bmatrix} b_0 \\ \vdots \\ b_m \end{bmatrix} \quad k \geq 0$$

$$c_j = 0, j < 0$$

Adding the $\breve{a}$ and $\hat{a}$ vectors of (3.28) and (3.29) gives the $a_{k+1} \cdots a_{k+m}$ coefficients as determined by (2.8). Combining these results with (3.27) reproduces the remaining portions of the equations in (2.8). The case for $k > m$ may be similarly demonstrated.

3.4 PM$^1$ Implementation, observations and performance

Successful implementation of this algorithm requires attention to numerics. While (3.14) stylistically uses the Moore-Penrose pseudoinverse, its well-known condition-number problems will cause this approach to underperform. A better numerical method is to perform QR factorization of the $C_2 = Q_2 R_2$ matrix, and then solve the eigenvalue problem as in (3.30), where the mathematical $R^{-1}$ operation is handled by performing backward substation of each column of $Q_2^T C_1$ through the upper-triangular $R$ matrix.

$$(3.30) \quad det(R_2^{-1} Q_2^T C_1 - \lambda I) = 0$$

Also, when solving (3.16), one only need solve the first $m$ equations using the first $m$ series coefficients. The result obtained solving the overdetermined problem using all $n = 2m + k + 1$ equations/coefficients is theoretically identical, but one does typically encounter more roundoff error. Finally, one does not need to find the residues at all. Instead, as demonstrated, one could construct the denominator polynomial and then solve for the numerator using the equations of (2.8).

The first set of validation tests involved duplicating the tests reported in section 2.2 using the PM$^1$ algorithm. The results are displayed in Table 3 for the $[m + k/m]$ $k \leq -1$ approximate of function $\hat{g}(z)$, with $2m + k$ series coefficients and compared with DM. The PM$^1$ algorithm shows no improvement over DM, nor should it as it should duplicate the results of DM, differing in only roundoff error. The algorithm with spurious pole information assimilation, PM$^2$ described in the next section, is able to remove all of the spurious poles for the 10 cases simulated. (These PM$^2$ results are reported here as it is convenient to compare the results in this table.)

Table 3 Poles for $[m + k/m], k \leq -1$ approximation of $\hat{g}(z)$

| Algorithm | System pole distance from 1 | Froissart Doublets near Unit Circle, $k \leq -1$ | Froissart Doublets near Unit Circle, $k > -1$ |
|---|---|---|---|
| DM | $O(\varepsilon)$ | $(m + k)$ | $(m - 1)$ |
| PM$^1$ | $O(\varepsilon)$ | $(m + k)$ | $(m - 1)$ |
| PM$^2$ | $O(\varepsilon)$ | 0 | 0 |

To test whether PM$^1$ results agree with the DM approach, we duplicated the simulations reported in section 2.2, using the noise model of (2.16) with uniformly distributed noise of level $\varepsilon$ on the interval $r_i \in [-1, 1]$ for:

- 10 cases of $[0/1]$ PAs with two different noise levels, $\varepsilon$, of uniformly distributed noise on the interval $r_i \in [-1, 1]$.
- 10 case of $[13/6]$ PAs with noise $\varepsilon = 10^{-6}$
- 10 case of $[5/14]$ PAs with noise $\varepsilon = 10^{-6}$

With the exception of variations in roundoff error, the PM$^1$ algorithm performed identically to PM.

Of the test we have used for validating PM$^1$, the most thorough one compares the performance of DM and PM$^1$ on the power series representation of the solution of the power-flow equations applied to the IEEE 118 bus-system [14], with modified loading selected to be at 91% of the voltage collapse point (saddle-node bifurcation point) to stress the algorithm. As a summative metric we use the maximum bus (nodal) power-flow mismatches taken over all buses for the near-diagonal $[m - 1/m]$ PA for $m \in \{1,\cdots,50\}, n = 2m$. The results of the maximum mismatches for DM and PM$^1$, plotted in Fig. 1, shows that performance is identical, except for small differences in roundoff error for $m>26$. The motivation behind selecting this problem and this metric is that for the bus maximum mismatches to agree, all PA's for the nodal/bus voltage variables for the 118 buses must show reasonable agreement. Any significant deviation would likely lead to a larger mismatch at some buses and this would likely show up as a deviation in the maximum mismatch.

The Maclaurin series for the voltage variables were calculated by holomorphically embedding the power flow equations using the canonical-form embedding with variable $\alpha$ as shown in (3.31)-(3.34), and then solving the linear recursion relationships for the series coefficients [15]. The traditional power-flow problem, models three sets of bus types, voltage magnitude controlled, no-voltage control, and one slack bus, denoted as $\{PV\}$, $\{PQ\}$ and $slack$ respectively. The definitions of the variables used are: $V_i(\alpha)$ is the voltage function at bus $i$, $P_i/Q_i/S_i$ is the real/reactive/complex-power injection at bus $i$, $Y_{ik}^{tr}$ are the elements of the admittance matrix with the shunt elements, $Y_i^{sh}$, removed, $V_{slack}$ is the slack-bus voltage and (*) is the complex-conjugate operator.

$$(3.31) \quad \sum_{k=1}^{N} Y_{ik}^{tr} V_k(\alpha) = \frac{\alpha S_i^*}{V_i^*(\alpha^*)} - \alpha Y_i^{sh} V_i(\alpha), \quad i \in \{PQ\}$$

$$(3.32) \quad \sum_{k=1}^{N} Y_{ik}^{tr} V_k(\alpha) = \frac{\alpha P_i - jQ_i(\alpha)}{V_i^*(\alpha^*)} - \alpha Y_i^{sh} V_i(\alpha), \quad i \in \{PV\}$$

$$(3.33) \quad V_i(\alpha) * V_i^*(\alpha^*) = 1 + \alpha(|V_i^{sp}|^2 - 1), \quad i \in \{PV\}$$

$$(3.34) \quad V_i(\alpha) = 1 + \alpha(V_{slack} - 1), \quad i \in slack$$

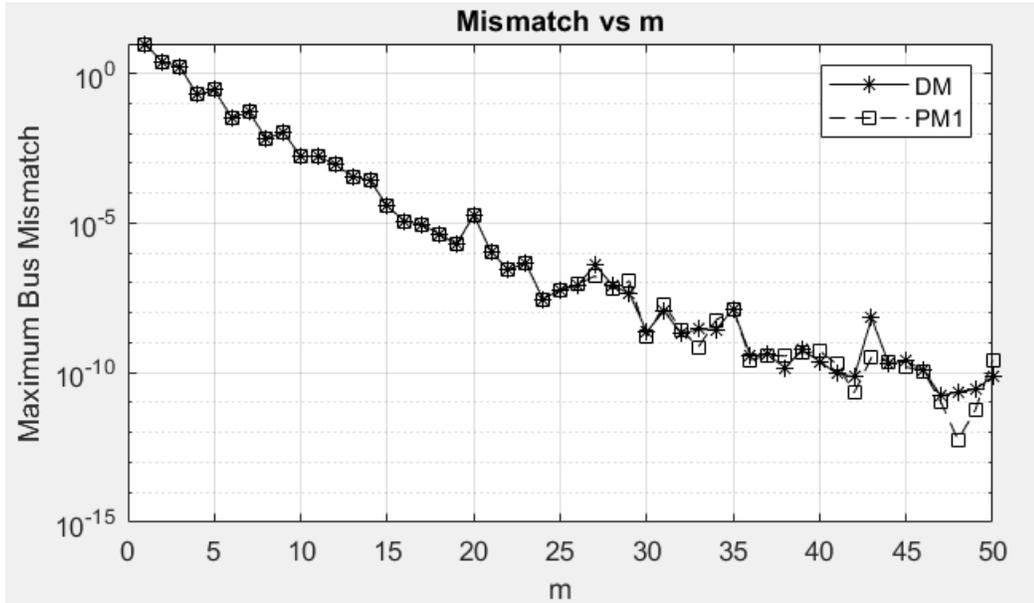

Fig. 1 Mismatch performance on the IEEE 118-bus system with modified loading.

## 3.5 PM[1] Limitations

Consider again the function in (2.4). The application of PM[1] in (3.35) produces the correct pole at zero. (This result also emphasizes why it is important to solve for the poles, rather than their inverse.)

(3.35)
$$C_1 = [c_1] = [0]$$
$$C_2 = [c_2] = [1]$$
$$det(C_2^+ C_1 - \lambda I) = det(0 - \lambda) = 0 \rightarrow \lambda = 0 = d_1^{-1}$$

Also, the correct numerator is obtained:

(3.36)
$$[c_1] = [1][e_1] \rightarrow [0] = [1][e_1] \rightarrow e_1 = 0$$

And finally using the equations of section 3.3,

(3.37)
$$[1/1] = c_o \frac{d_1^{-1} - z}{d_1^{-1} - z} + \frac{d_1^{-1} e_1}{d_1^{-1} - z} = \frac{-z}{-z} + \frac{0}{0 - z} = \frac{z}{z}$$

This example belies the fact that, like DM, PM[1] suffers from a degeneracy problem: if a pole occurs at/near zero in the case where $m$ is greater than 1, then the matrix of (3.16), or (3.26), will/can become numerically singular. Even for poles not terribly close to zero, if $m$ becomes large enough and/or if one solves the overdetermined problem, the $D$ matrix in (3.16), or (3.26), can become numerically rank deficient.

Given that DM is more numerically efficient than PM[1], we are hard pressed to see the advantages for this more complicated formulation; however, as we will see in the next section, a modification of this approach which will lead to the ability to identify spurious poles caused by precision limitation, poles place at/near the origin due to unfortunate selection of $m$ or $k$ and a means of adjusting $m$ or $k$ to eliminate these poles, including spurious poles inside the radius of convergence, without eliminating the information carried in the series terms responsible for these poles.

## 4 Padé Matrix Pencil Method with Spurious Pole Information Assimilation: Noisy Data with Precision Limitation (PM[2] aka PM2)

We have observed when working on the power-flow problem for electric power transmission that, despite the onset of spurious poles as the number of series terms increases, the accuracy of the functions under study often continue to improve, indicating that the spurious poles contain information about the underlying function, even though they may not lie on the branch-cut for this problem. Our goal has been to find an algorithm that could assimilate this information without producing so-called spurious poles. In fact, because some of these poles contain information about the underlying function, calling them spurious may be a misnomer to some degree.

In building a PA, if our selection of $m$ (or $k$) leads to a defect, $\gamma_{k,m} > 0$, this defect is typically corrected by reducing $m$ (or modifying $k$). This may be accomplished by either reducing $n$, effectively discarding the information of the rejected coefficients, or keeping $n$ unchanged and reducing $l$ below $m$, so that the number of poles becomes $l$, as discussed below.

Provided our PA has a zero defect, the poles of the PA of conformation, $[m + k/m], -m \leq k \leq n - 2m - 1$, are calculated using the generalized eigenvalue problem of (3.14) and residues may be found using (3.16) or (3.26) depending on the conformation. The $C_1$ and $C_2$ matrices are defined by (3.22), which are submatrices of the $C$ matrix in (3.20). In the previous sections, while the theoretical derivation was more general, we considered only the cases where we enforced the equality constraint, $n - l = l = m$. In assimilating, rather than rejecting, the information provided by series coefficients that are responsible, in a sense, for spurious pole production, it is helpful to develop an intuitive understanding of the implication of relaxing this equality constraint, $n - l \neq l \neq m$. Relaxing this equality constraint while restructuring the eigenvalue problem, will allow us to develop an algorithm with the capabilities that we need.

It is insightful to consider the implications of the relaxed constraint of two distinct problem types: functions whose exact representation involves a finite and then an infinite number of poles. We will consider each of these in turn, assuming precision is not a limitation. The insight gained will inform how coefficient perturbations, due to measurement noise or finite precision calculations may be accommodated. Because of page limitations, we consider in detail only the $[m − 1/m]$ case. The cases for other conformations may be described similarly.

4.1  Finite number of poles

Assume that we are dealing with a meromorphic function with a finite number of poles, $m$, exactly expressed with a PA of conformation $[m − 1/m]$. Using the previous (or any number of) PA-producing algorithm with $n = 2m$ series coefficient, we obtain the theoretically exact PA as shown in Table 4, case (a) with set of poles $p = \{p_1, p_2, \cdots p_m\}$.

Assuming we do not know the true value of the number of poles a priori, and instead select $m', n' = 2m'$. If $m' < m$ and we select $n − l = l = m'$, Table 4, case (b), then the previous algorithm will perform equivalently to other algorithms, generating a set of poles corresponding to the equivalent denominator polynomials with poles $\{\hat{p}_1, \hat{p}_2, \cdots \hat{p}_{m'}\} \neq p$. If we select $l < m'$, Table 4, case (c), then, given the dimension of $[C_2^+ C_1]_{l \times l}$, the number of eigenvalues and poles calculated will be $l < m'$ and the degree of the resultant denominator polynomial will be reduced below our selected $m'$. In fact, that is the motivation for the restriction of $m \leq l \leq n − m$ in section 3.1, a restriction which is only meaningful when we are without defect.

If we select $m'' > m$, $n'' \geq 2m''$, $m'' \geq l > m$, given the dimension of $[C_2^+ C_1]_{l \times l}$, the algorithm of the previous sections will produce the desired $m$ poles and $l − m$ poles at the origin, as shown in Table 4, case (d). This is one example of an improper conformation selection that results in spurious poles at the origin, if exact arithmetic is assumed, or near the origin if roundoff error is encountered. This situation is particularly informative because it suggests two methods for dealing with these superfluous poles: either the number of series coefficients can be reduced from $2m''$ to $2m$, discarding the information these coefficients contain, or the value of $l$ reduced to $m$, assuming one knows $m$ or can estimate it. (We will demonstrate later how $m$ may be estimated.) Assuming exact arithmetic and that the superfluous coefficients are consistent with the definition of the function, not erroneous, reducing $l$ from $m''$ to $m$ effectively reconfigures the algorithm to have $m'' − l = m'' − m$ redundant relationships, eliminating the $\gamma_{-1,m} = 2(m'' − m)$ defect, with no degradation of the algorithm's performance, numerics notwithstanding.

Table 4 Poles for $[m − 1/m]_{f(z)}$

| Case | (a) | (b) | (c) | (d) |
|---|---|---|---|---|
| $m$ | $m$ | $m' < m$ | $m' < m$ | $m'' > m$ |
| $n$ | $2m$ | $2m'$ | $2m'$ | $n \geq 2m''$ |
| $l$ | $l = m$ | $l = m'$ | $l < m'$ | $m'' \geq l > m, l < n$ |
| Poles | $p = \{p_1, p_2, \cdots p_m\}$ | $\{\hat{p}_1, \hat{p}_2, \cdots \hat{p}_{m'}\}$ | $\{\tilde{p}_1, \tilde{p}_2, \cdots \tilde{p}_l\}$ | $\{p_1, p_2, \cdots p_m, 0, \cdots, 0\}_{1 \times l}$ |

4.2  Infinite number of poles

Table 5 Poles for $[m − 1/m]_{f(z)}$, infinite number of poles

| Case | (a) | (b) |
|---|---|---|
| $m$ | $m' < \infty$ | $m' < \infty$ |
| $n$ | $2m'$ | $n \geq 2m'$ |
| $l$ | $m'$ | $m' \geq l$ |
| Poles | $\{\hat{p}_1, \hat{p}_2, \cdots \hat{p}_{m'}\}$ | $\{\tilde{p}_1, \tilde{p}_2, \cdots \tilde{p}_l\}$ |

If we are dealing with a meromorphic function with an infinite number of poles, or a function with branch points, then regardless of the finite value of $m'$ and/or $l$ we choose, the set of $m$ poles produced by the algorithm of the previous section will always be insufficient to exactly represent the function, but will serve as an approximation as indicated in Table 5, whose set of spurious poles will be dependent on the PA conformation chosen. However, as we increase the number series terms, the increasing number of floating-point calculations typically involved in calculating the higher order terms introduces roundoff error, which appears as a noise perturbation in the series coefficients, leading to spurious poles occurring in the PA regardless of the conformation chosen.

### 4.3  Spurious pole recognition and information assimilation

The $l$ parameter can be used to filter out spurious poles by adjusting the assumed $m$ value, which changes the conformation of the PA. However, the first step is recognizing the existence of spurious poles and being able to estimate the degree of the defect, due to either noise in the series coefficients or improper conformation.

Consider the $C$ matrix constructed using the definitions in (3.8), generalized for an arbitrary $[m + k/m]$ conformation, as in (3.20) and then further generalized to allow for the case where the number of poles has been reduced to $l$ to eliminate the current estimate of the defect of $\gamma_{k,m} = 2(m - l)$. This changes the conformation of the C matrix from $(m) \times (m + 1)$ to $(2m - l) \times (l + 1)$. The SVD of this C matrix is given by,

$$(4.1) \quad C_{2m-l \times l+1} = U_{2m-l \times 2m-l} \Sigma_{2m-l \times l+1} V_{l+1 \times l+1}^H$$

where $U$ and $V$ are unitary and $\Sigma$ is a real-valued diagonal matrix with singular values $\sigma_1 \geq \sigma_2 \geq \cdots \sigma_{l+1} \geq 0$, $2m - l > l + 1$. Using exact arithmetic, if the defect is greater than $\gamma_{k,m} = 2(m - l)$, say $l' < l$ then all of the eigenvalues indexed, $l', l' + 1, \cdots, l$, will be zero. However, when noise is added to the series coefficients, these values will be greater than zero. We have observed that the singular values leading to spurious poles may be identified by their magnitude. Specifically, let $\sigma_{max}$ be the maximum singular value of $\Sigma$ and let $t$ be the number of accurate digits expected in the series coefficients. We have observed that when the $\Sigma$ matrix is truncated to eliminate all singular values, $\sigma_i$ such that,

$$(4.2) \quad \frac{\sigma_i}{\sigma_{max}} < 10^{-t} \quad t > 0$$

then the spurious poles due to those singular values are eliminated. (This is an interesting result because it suggests that the noise level contained in any series may be estimated by adjusting $t$ to the minimum value that eliminates all spurious poles due to noise.) This requires adjusting the conformation of the $C_1$ and $C_2$ matrices with the resultant conformation of the $C$ matrix becoming $C_{2m-l' \times l'+1}$ and a new SVD for this conformation must be performed, iteratively using the test of (4.2) until no further singular values exceed this criterion.

When this iterative process has been completed, we now have the correct conformation to apply to the $C_1$ and $C_2$ matrices in (3.8) to eliminate the spurious poles due to noise in the series coefficients for the initial $[m + k/m]$ conformation selected and the eigenvalues may be found using (3.30); however, we can use the final SVD in this iterative process to eliminate some of the complexity of the QR factorization step needed in (3.30).

If we eliminate the first (last) column of $V^H$ in (4.1), the resultant matrix multiplication in (4.1) correspond to the $C_1$ ($C_2$) matrix. Letting $V_1^H(V_2^H)$ corresponds to the first (last) $l$ columns of $V^H$ the $C_1$ ($C_2$) matrix may be written.

$$(4.3) \quad C_{1(2m-l \times l)} = U_{(2m-l \times 2m-l)} \Sigma_{(2m-l \times l+1)} V_{1(l+1 \times l)}^H$$

$$(4.4) \quad C_{2(2m-l \times l)} = U_{(2m-l \times 2m-l)} \Sigma_{(2m-l \times l+1)} V_{2(l+1 \times l)}^H$$

Observe that we can arrive at the equivalent eigenvalue problems as follows.

$$(4.5) \quad \begin{aligned} (C_1 - \lambda C_2)_{2m-l \times l} &\to (U\Sigma)_{2m-l \times l+1} \left(V_1^H - \lambda V_2^H\right)_{l+1 \times l} \to \left(V_1^H - \lambda V_2^H\right)_{l+1 \times l} \\ &\to \left(V_2^{H+} V_1^H\right)_{l \times l} - \lambda I \qquad l \leq m \end{aligned}$$

This is the implementation we have used in the work reported here. This saves some computation time since the dimension of $V_{2(l+1 \times l)}^H$ and $V_{1(l+1 \times l)}^H$ is less than that of $C_{2(2m-l \times l)}$ and $C_{1(2m-l \times l)}$ and QR factorization of $V_2^H$ along with solution of the overdetermined problem $R_{V2}^{-1} Q_{V2}^T V_1^H$ in (4.6), will incur less roundoff error than for the larger matrix, $C_1$; however, roundoff error is incurred by arriving at $V_2^H$ and $V_1^H$. We have not studied the trade off in these two approaches.

$$(4.6) \quad det\left(R_{V2}^{-1} Q_{V2}^T V_1^H - \lambda I\right) = 0$$

Once the eigenvalues are calculated using in (4.6), the eigenvalues may be inspected to find the ones placed near the point of development (zero in our case) that are spurious, either because of precision issues or because of an unfortunate selection of $(k, m)$. These poles may be removed using the following procedure. (Presumably the user has scaled the series so that the ROC is near 1.0 but, regardless, if the ROC is known, the following procedure can also be used to remove any poles within the ROC.)

Because the theory developed only strictly applies when we keep the Hankel matrix structure, if a pole near the origin (or spurious pole within the radius of convergence) is to be removed, this is accomplished by changing the conformation of the the $C_1$ and $C_2$ matrices, using the same procedure for dealing with noise in the series coefficients: adjusting the $l$ parameter to effectively reduce the number of columns by 1 and increase the number of rows by 1, and then performing an SVD on the resultant matrix to ensure that no spurious poles due to precision have been introduced. (Though our experience is limited, we have not seen precision-induced spurious poles introduced during this step but, as described later, depending on the implementation, it is possible for this to occur.) If noise-induced spurious poles remain or are created by this step, i.e., singular values that violate (4.2) are encountered, revision of the conformation of the $C_1$ and $C_2$ matrices is continued until no singular values violate the (4.2) criterion and the eigenvalues/poles of the PA with the corresponding revised conformation are again calculated. This scheme is repeated until none of the poles we wish to remove by this procedure remain. It is imperative that all poles near zero be eliminated, otherwise the $D$ matrix in (3.26) can become numerically singular.

### 4.4  Residue calculation

Once the undesired poles have been removed, (3.16) or (3.26) is solved for the residues, $-e_j d_j^{-1}$. Assuming there has been an adjustment to number of poles, from the ideal, $m$, to a number without spurious poles, $l < m$, the entire overdetermined problem is solved for the residues, unlike in the ideal case where only the first $m$ coefficients need be used.

Because the coefficients of the $D$ matrix grow/shrink exponentially with the $n$, if small poles are not removed, (e.g., in the desire to retain spurious poles unrelated to noise or the unfortunate selection of $(k, m)$, i.e., the spurious poles mentioned in section 2.3) the $D$ matrix can become numerically singular. To that end, before constructing the $D$ matrix, we discard any eigenvalues corresponding to poles within a circle whose radius is $10^{-3}$, centered at the origin. Because of the possibility of the $D$ matrix becoming numerically rank deficient, in our implementation we perform an SVD of the $D$ matrix to identify any small singular values, using the criterion of (4.2) to determine if the conformation of the $C$ matrix should be changed. If small singular values are encountered, the number or columns/rows of the $C$ matrix is reduced/increased by one. We do not eliminate any poles beyond the ROC. If this process fails to keep the $D$ matrix from becoming numerically singular (a situation which we have not encountered) the number of equations/coefficients used in $c = De$ may be reduced until the $D$ matrix becomes numerically nonsingular. We have found this last step to be unnecessary if all poles inside the ROC are removed.

### 4.5  The PM$^2$ algorithm

If the noise level results in the number of accurate digits being significantly less than the precision level of the computing engine, then the roundoff encounter by the mathematics presented in sections 4.3-4.4 occupies digits far from the accurate digits of the series coefficients and the algorithm tends to behave as if it is using exact arithmetic. If the number of accurate digits approaches that of machine precision, then the interaction of roundoff error in the series coefficients and roundoff error produced by the algorithm interact and the outcome becomes more sensitive to the implementation. The most demanding task is the removal of spurious poles due to precision limitation because the roundoff perturbation in the series coefficients spans a few orders of magnitude that is shared by the roundoff error of the algorithm.

In this section we present the algorithm capable of removing spurious poles due to precision limits. In the next section, we address several short cuts along with comments about their anticipate impact on the accuracy of the results.

The goal of removal of spurious poles is, frankly, an ill-posed goal. There may be many ways of changing the confirmation of the PA which will eliminate spurious poles with more or less (or essentially the same) accuracy. For example, if an $[m + k/m]$ PA is found to have no spurious poles due to noise but has one pole at the origin due to an unfortunate selection of $(k, m)$, the algorithm presented below reduces $m$ by one, while keeping $k$ fixed, and continues reducing the denominator polynomial degree until the pole near the origin is gone or the degree reaches zero. Recognizing that the $C$ matrix always requires an even number of coefficients, one could also change the conformation by decreasing $k$ by 2 and increasing $m$ by 1 to eliminate the pole at the origin. Making informed decisions in this regard requires knowledge of the Padé table [1], [16], [17], which is rare to have.

In the algorithm below we arbitrarily assume that $k$ is to be preserved, which may prevent the elimination of all spurious poles in some cases. It is relatively easy to adjust the algorithm to allow $k$ to change and preserve $m$, or to allow both to change. The user will want to modify how the algorithm adjusts the conformation of $C_1$ and $C_2$ matrices to suit their needs.

**PM² Algorithm**

Input: $m \geq 0, k \geq -m, 2m + k = n$, series coefficients $c_0, \cdots, c_{2m+k}$, each with $t > 0$ accurate digits of precision, with series scaled so that the ROC is near 1.

Output: Denominator coefficients, $b_0, \cdots, b_l$ and numerator coefficients, $a_0, \cdots, a_{k+l}$

1. For a given selected $[m + k/m]$ conformation, set $l = m$
2. Construct $C_{(l \times l+1)} = C_{(2m-l \times l+1)}$
3. Calculate the SVD of $C_{(2m-l \times l+1)}$
4. If $l = 1$, jump to Step 7.
5. Count the number of singular values, $n_s$ such that $\frac{\sigma_i}{\sigma_{max}} < 10^{-t}$.
6. If $n_s > 0$, set $l = l - n_s$ return to Step 2. Otherwise, proceed.
7. Calculate the eigenvalues from $det(R_{V2}^{-1} Q_{V2}^T V_1^H - \lambda I) = 0$. If any eigenvalues are zero or within a user defined radius, e.g., 10⁻³, remove these eigenvalues before building the $D$ matrix.
8. Calculate the singular values of $D$. If the smallest singular value is less than $10^{-t}$ times the maximum singular value, set $l = l - 1$ return to Step 2. Otherwise proceed.
9. Calculate the residues using $c = De$, (3.16) or (3.26).
10. From the poles, residues and series coefficients, calculate the numerator and denominator polynomial coefficients if necessary.

4.6    Short cuts

In cases where the expected number of accurate digits in the Maclaurin series coefficients is much below machine precision and/or the number of series coefficients is not too great, the algorithm as presented in the previous section may be simplified. Below is a list of computational-complexity shortcuts that work reasonably well in many cases. We have limited experience with these short cuts but share that experience below.

- When calculating the residues, use a number of coefficients less than the total number of series coefficients but at least equal to the degree of the denominator. Alternatively, calculate the coefficients of the denominator polynomial and use (2.8)(b) or (2.8)(c). We have seen both of these approaches work well provided the number of spurious poles is not large.
- When a pole (from solving the eigenvalue problem of (4.5)) is sufficiently near the origin, rather than change the conformation of the $C$ matrix, simply eliminate the pole from the $D$ matrix and solve for the residues, recognizing that denominator degree has been reduced by one. This is an approximation because the relation implied by the Hankel formulation is violated.
- The following does the most violence to the method and does not eliminate all of the poles due to noise, nor does it necessarily deal with any poles near the origin, but it is much faster and does provide one mechanism for eliminating coefficients whose information does not improve the accuracy of the result. In this approach, the singular value criterion of (4.2) is performed once. Assume the number of singular values (poles) has been reduced from $m$ to $l$. The $V_1'^H (V_2'^H)$ matrices corresponding to the first (last) $l$ columns of $V'^H$. Note that construction of $C'_1$ using $V_1'^H$ will no longer yield a Hankel matrix. A slightly improved method is to use $l$ to inform the conformation of the $C_1$ and $C_2$ matrices and use (3.30) to find the eigenvalues/poles and then adjust the conformation if poles near the origin are found.

## 5 PM² Numerical Results

### 5.1 Geometric Series

Obtaining closed form estimates of the behavior of PA algorithms under the influence of series-coefficient noise is usually not practical. However, this is possible for the special case of a PA with conformation [0/1] when PM² is applied to a geometric series corrupted by the noise described in (2.16). It is possible to show (not developed here) that the system pole (related to the analytic properties of the function), numerator coefficient, and function value errors are $O(\varepsilon)$, assuming $r$ is uniformly distributed on the interval $r_i \in [-1, 1]$.

To verify the pole error, shown in Fig. 2 are plots of the average of the magnitudes of the system pole error averaged over a sample size of 10 for both the DM and PM² algorithm using a series of $n = 20$ with the noise model described in (2.16). The DM algorithm produces a [9/10] PAs from each series, which uses a noise amplitude of $\varepsilon$. The PM² algorithm builds a [0/1] PA after spurious pole information assimilation with parameter $10^{-t} = \varepsilon$. This plot shows the superior pole error performance that can be achieved when superfluous noisy series coefficients are used as redundant data.

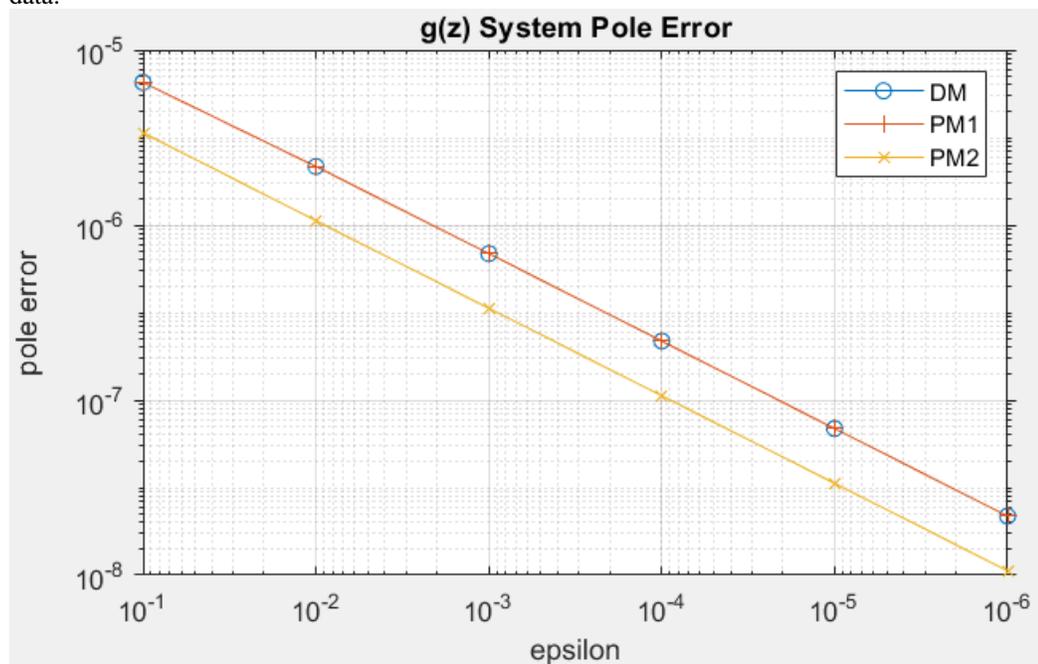

Fig. 2 $g(z)$ PA system pole error versus noise level for $n=20$, with a sample size of 10 using DM, PM¹ and PM², $10^{-t} = \varepsilon$.

Using the same data set ($n = 20$, sample size of 10) the experimental function error in the $x$-range $[-0.9, 0.9]$, $[0.9, 0.99]$, $[1.1, 100]$ is shown in Fig. 3-Fig. 5 with noise level as a parameter for both the DM and PM² methods, again with $10^{-t} = \varepsilon$. These figures show that the PM² function error, using a [0/1] PA, is at least as accurate, and often times superior to DM, a [9/10] PA. This plot also shows that the PM² error is immune to spikes caused by spurious poles. We have observed similar results with more complex meromorphic function.

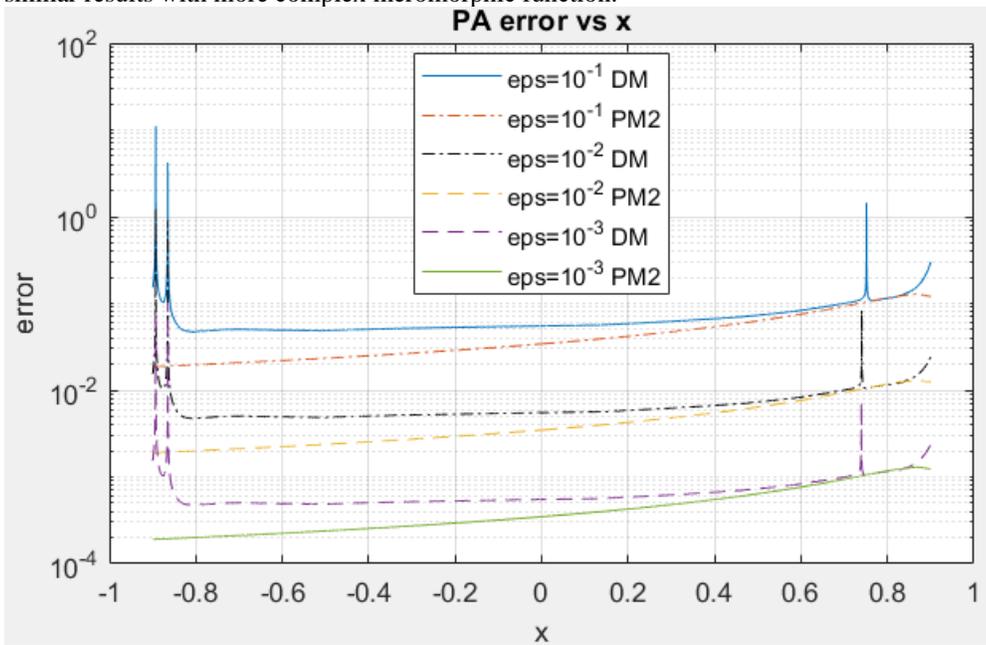

Fig. 3 $g(z)$ PA error versus $x \in [-0.9, 0.9]$, $y=0$ with noise level as a parameter using DM and PM², $10^{-t} = \varepsilon$.

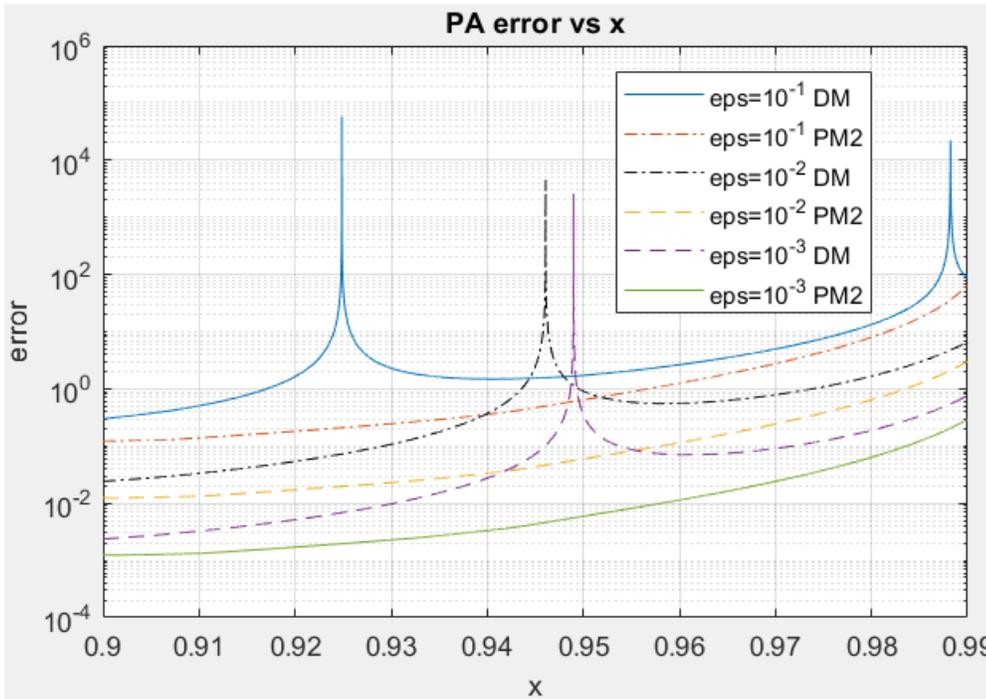

Fig. 4 $g(z)$ PA error versus $x \in [0.9, 0.99]$, $y=0$ with noise level as a parameter using DM and PM², $10^{-t} = \varepsilon$.

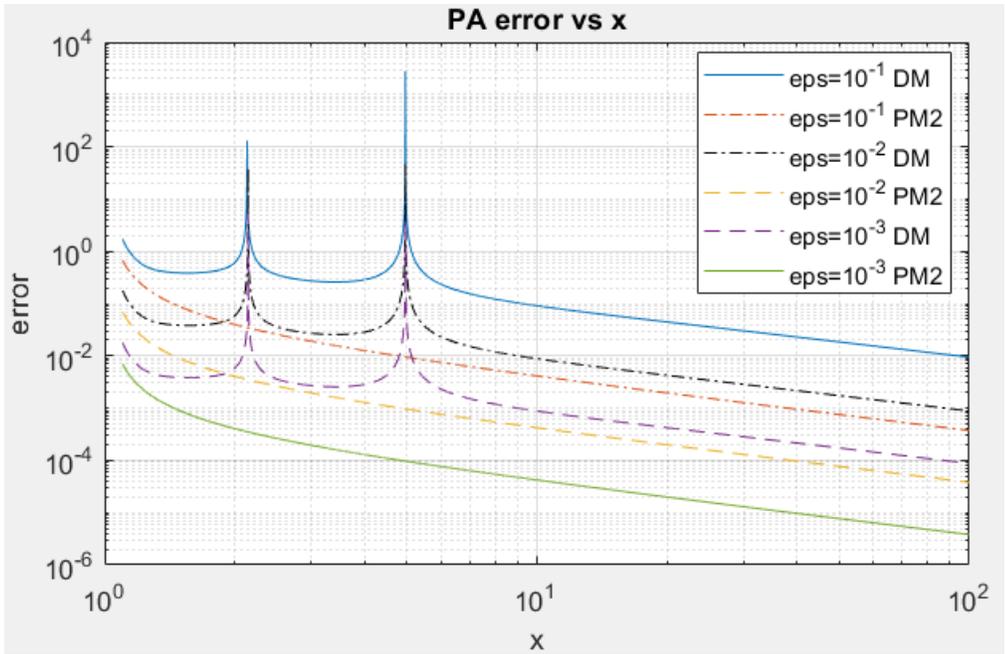

Fig. 5 $g(z)$ PA error versus $x \in [1.1, \ 100]$, $y=0$ with noise level as a parameter using DM and PM², $10^{-t} = \varepsilon$.

## 5.2 $Log(1.2 - z)$

Fig. 6 shows the PA poles resulting from using an $n=41$ series representation of the function $f(z) = log(1.2 - z)$, using both the DM algorithm, which produced a [20/20] with 6 spurious poles, and the PM² algorithm, with $t=14$, which produced an [12/12] PA with no spurious poles. The maximum error of the PAs in the unit circle discretized into uniform orthogonal meshes where the distance between adjacent mesh points is 0.01, and one mesh point anchored at the origin (7860 mesh points), was $2.9 \times 10^{-12}$ for DM and $1.6 \times 10^{-11}$ for PM2. In Fig. 6 and all other pole-zero plots in this work, we use the convention that the marker "x" is for poles and the marker "o" is for zeros. Note that according to Stahl's theorem for functions with branch points, the poles characteristic of the underlying function, accumulate on a branch cut which has the specific property of minimal logarithmic capacity, a branch cut that connects subsets of the function's branch points [4]. For this $log$ function, the functions poles lie along this branch cut, which is a ray coincident with the real line, starting at $x = 1.2$.

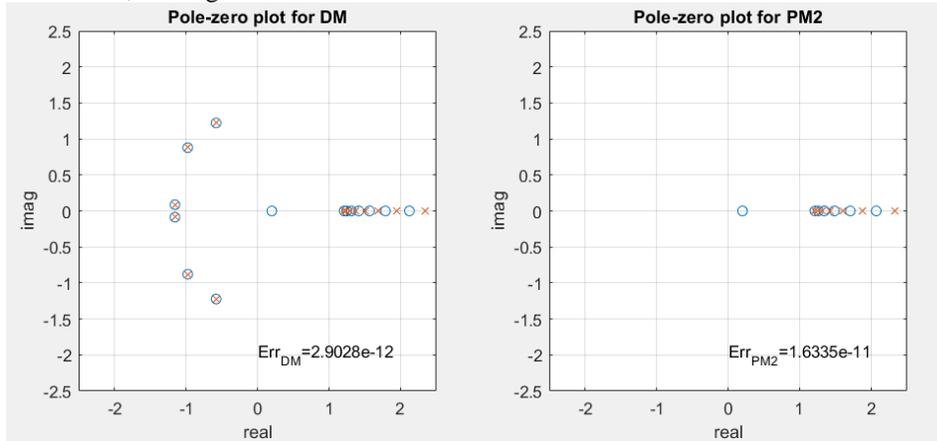

Fig. 6 Poles and zeros of $log(1.2 - z)$ using DM and PM² with a value of $t=14, n = 40$.

The accuracy with which a function is rendered may improve as more terms are added to the series that characterizes the function, even if spurious poles are added by these terms. An important question to ask is: In this case, do the added spurious poles contribute accuracy to the function, or is the enhanced accuracy caused by the adjustment to the location of those poles and zeros (related to the analytical properties of the function) caused by the added series terms? Shown in Fig. 7 is the PA error plot for $log(1.2 - z)$ for a range along the real axis of $x \in [0, 1]$, with $n = 40$, using PM$^1$ (all poles), PM$^1$ with spurious poles removed and PM$^2$(with spurious pole information assimilated, $t = 14$). This plot demonstrates that, in some cases, the spurious poles contribute to the accuracy of the function and that PM$^2$ assimilates some of this information into the reduced order PA.

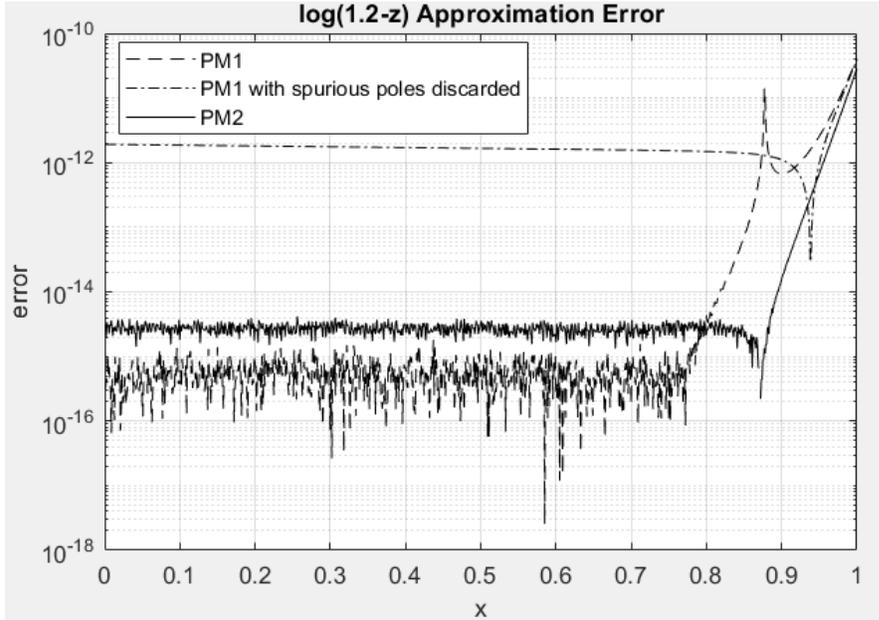

Fig. 7 $log(1.2 - z)$ error versus $x \in [0, 1]$, y=0, using PM1, PM1 with spurious poles discarded and PM2 (with spurious pole information assimilation), t=14, $n = 40$.

5.3  The 118-bus power-flow problem

The 118-bus [14] power flow problem defined by (3.31)-(3.34), loaded at 91% of the saddle-node bifurcation point, was solved using DM, RPA (with $tol = 10^{-14}$) [6], [7], and PM$^2$ with $t$=14. As a summative metric we use the bus (nodal) maximum mismatches at this loading level, taken over all buses for the near-diagonal $[m - 1/m]$ PA for $m \in \{1, \cdots, 50\}, n = 2m$. The results of the maximum mismatches and number of poles retained (for representative bus 22 only) for DM, RPA and PM$^2$ is shown in Fig. 8 versus $m$. This plot reveals several phenotypical features of the algorithm. Most notably, the accuracy of the PM$^2$ solution is typically as good as DM. Also, the progress toward an accurate solution tends to be more well behaved than the progress of DM. The number of retained poles is a monotonic function, with a maximum value of 26 for this experiment, whereas the equivalent maximum number for DM is 50. The fact that the accuracy increases as we add more series terms but do not add more system poles is consistent with the claim that we are assimilating into the PA information that would be contained in the spurious poles. This figure shows that, on this problem, PM$^2$ out-performs RPA in terms of accuracy and minimizing the number of poles need to represent the voltage functions, while also being capable of representing the 117 complex-valued voltage functions with the accuracy similar to that of DM, where accuracy is measured using bus mismatches. The pole-zero plots of the PAs obtained using DM, RPA and PM$^2$, in Fig. 9 for bus 22 (which is representative of the other buses), shows that PM$^2$ is able to eliminate the spurious poles within the ROC using a value of $t$=14, indicating that the recursion relationships used to generate the $V(\alpha)$ series preserved at least 14 of the 15.7 digits of accuracy used in the calculation. This figure also shows that PM2 and RPA give similar results, removing spurious poles produced using DM.

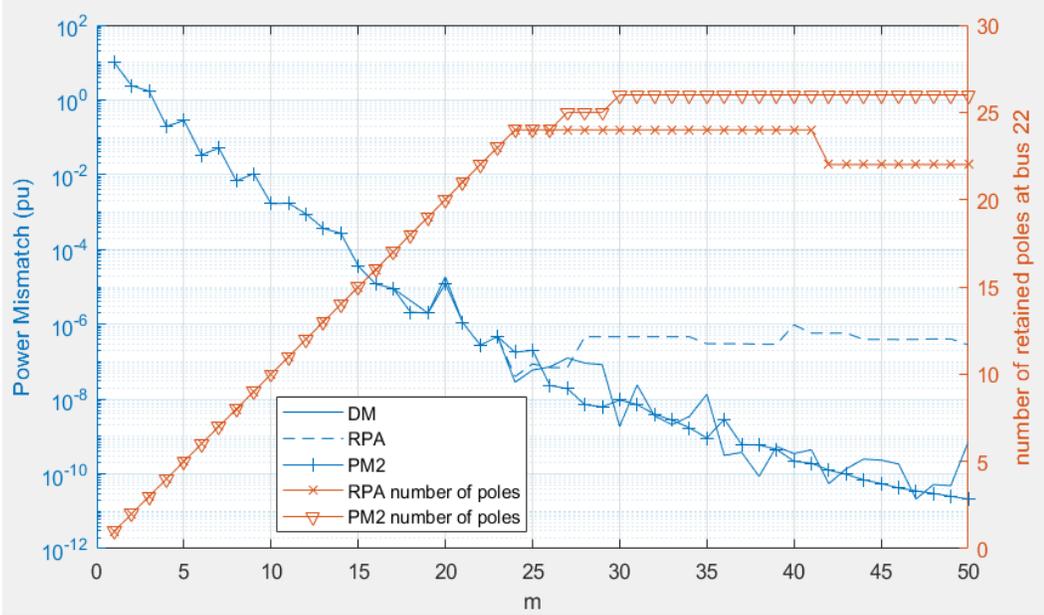
Fig. 8 Nodal power mismatches and number of retained poles versus $m$ for the 118-bus power flow.

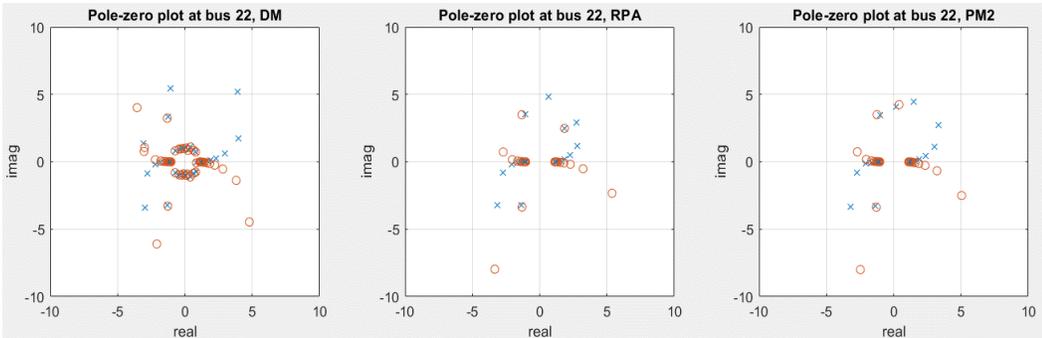
Fig. 9 PA Poles and zeros of bus 22 for the DM, RPA and PM$^2$ algorithms

## 6    Conclusion

We have developed the theory of a new method. PM[1], for calculating Padé approximants that is capable of matching the results of the well established $O(m^3)$ direct method (cf. [1], Ch. 2), including producing the same set of spurious poles; any numerical differences between the method are due to roundoff error. A generalization of the approach, PM$^2$, aids in the elimination of spurious poles inside the radius of convergence, including those responsible for degeneracy in the direct method. The approach can also be used to eliminate spurious poles due to precision limitations or coefficient noise. The method is more computationally complex than the direct method, but allows the information about the underlying function contained in the spurious poles to be assimilated, resulting in a reduced order PA. Numerical experiments show it to be competitive with the direct method in terms of accuracy for the cases presented while significantly reducing the denominator and numerator degrees. The approach is also shown to be competitive with the performance of the Robust Padé Approximation method.

The one weakness of the algorithm as presented is that no data are discarded. This means that if an exorbitant number of series terms are generated, resulting in gross errors due to precision limitation, or if the noise level of the series coefficients is so high as to produce occasional bad data point, these terms used redundantly can skew the estimates. Methods for discarding data are easily envisioned, which use a modification to the filter parameter portion of the algorithm, but are not presented here.

# 7 Acknowledgement

We wish to thank Professor Emeritus Ronald Roedel and Michael Thomas Masengarb for helpful discussions during the execution of the research. Much obliged. A special thanks to Abhinav Dronamraju who did some initial work on degeneracy before graduation and Qirui Li whose work recognizing patterns in condition number behavior was pivotal.